\documentclass[12pt]{article}

\usepackage{amssymb}
\usepackage{amsmath}
\usepackage[latin1]{inputenc}
\usepackage[dvips]{epsfig}

\newtheorem{thm}{Theorem}[section]
\newtheorem{prop}[thm]{Proposition}
\newtheorem{lem}[thm]{Lemma}

\newtheorem{cor}[thm]{Corollary}

\newcommand{\g}{{\gamma}}

\newcommand{\be}{{\beta}}

\newcommand{\s}{{\sigma}}

\newcommand{\SO}{\operatorname{S}}
\newcommand{\sign}{\operatorname{sign}}

\newcommand{\Neg}{\operatorname{Neg}}

\newcommand{\maj}{\operatorname{maj}}
\newcommand{\fmaj}{\operatorname{fmaj}}
\newcommand{\Dmaj}{\operatorname{Dmaj}}
\newcommand{\inv}{\operatorname{inv}}
\newcommand{\me}{\operatorname{N}}

\newcommand{\des}{\operatorname{des}}
\newcommand{\Des}{\operatorname{Des}}

\newcommand{\NN}{\mathbb{N}}
\newcommand{\ZZ}{\mathbb{Z}}
\newcommand{\PP}{\mathbb{P}}

\newenvironment{proof}{\begin{trivlist}\item{\bf{Proof.}}}
  {\hfill\rule{2mm}{2mm}\end{trivlist}}

\begin{document}

\title{Signed Mahonian polynomials \\ for classical Weyl groups}
\author{Riccardo Biagioli}

\date{\today}
\maketitle

\begin{abstract}
The generating functions of the major index and of the flag-major index, with each of the one-dimensional characters over the symmetric and hyperoctahedral group, respectively, have simple product formulas. In this paper, we give a factorial-type formula for the generating function of the $D$-major index with sign over the Weyl groups of type $D$. This completes a picture which is now known for all the classical Weyl groups.
\end{abstract}

\section{Introduction}
Sums of the form 
\[\sum_{w\in W} \chi(w)q^{\des(w)},\]
where $W$ is a classical Weyl group, $\chi$ is a one-dimensional character of $W$, and $\des(w)$ is the number of descents of $w$ as Coxeter group element, have been investigated by Reiner \cite{Re}. In the case of the symmetric group, when $\chi$ is the trivial character this sum is the well known {\em Eulerian polynomial} \cite{Ca}, and when $\chi$ is the sign character then it is the {\em signed Eulerian polynomial} studied by D\'esarm\'enien and Foata \cite{DF}, and Wachs \cite{W}. Analogously, consider the sum
\begin{equation}\label{character}
\sum_{w \in W} \chi(w) q^{\maj_W(w)},
\end{equation}
where $\maj_W$ denotes a suitable {\em Mahonian statistic} on the corresponding group $W$. Recall that a statistic on a Coxeter group is said to be Mahonian if it is equidistributed with the length function on the group. In the case of the symmetric group, if $\chi$ is the trivial character, the sum in (\ref{character}) is nothing but the Poincar\'e polynomial of $S_n$, which, as is well known, admits a nice product formula for every finite Coxeter group (see e.g., \cite{Hum}). Otherwise, if $\chi$ is the sign character, this sum corresponds to the  {\em signed Mahonian polynomial} studied by Gessel and Simion \cite{W}, who found an elegant product formula for it in terms of $q$-factorials. Recently, several extensions of this result have been given by Adin, Gessel and Roichman \cite{AGR}. In particular they have provided nice formulas for the polynomial in (\ref{character}) in the case of the hyperoctahedral group $B_n$ equipped with the Mahonian statistic ``$\fmaj$'', the {\em flag-major index}, which was defined by Adin and Roichman in \cite{AR}.

In this paper we deal with Weyl groups of type $D$ together with the {\em $D$-major index} ``$\Dmaj$'', defined by the present author and Caselli in \cite{BC}. The $D$-major index is a Mahonian statistic that has the analogous role for $D_n$, as $\maj$ has for $S_n$ and $\fmaj$ for $B_n$. Moreover it shares with them very nice algebraic properties (see \cite{GG},\cite{AR},\cite{ABR,ABR1} and \cite{BC,BC2}). Like the symmetric group, $D_n$ has only two one-dimensional characters, the trivial and the sign. In the case of the trivial character the corresponding sum of (\ref{character}) is again the Poincar\'e polynomial of $D_n$. In the case of the sign character we give a factorial-type formula for the signed Mahonian polynomial of type $D$, extending the formulas previously mentioned. Toward this end, we define a natural sign-reversing involution on $B_n$, in the ``style of Wachs'' (see \cite{W}), that is $\fmaj$ preserving. We use this involution, first to give an easy proof of the Adin-Gessel-Roichman formula for the signed Mahonian polynomial for $B_{2n}$, and then to derive from the latter formula the analogue for $D_n$. 
This completes a picture for the generating functions for the major index with one-dimensional characters over the classical Weyl groups $S_n$, $B_n$, and $D_n$.
  
\section{Preliminaries and Notation}

In this section we give some definitions, notation and results that
will be used in the rest of this work. For $n \in \NN$ we let $[n]:= \{ 1,2, \ldots , n \} $ (where $[0]:= \emptyset $). Given $n, m \in \ZZ, \; n \leq m$, we let $[n,m]:=\{n,n+1, \ldots, m \}.$ We let $\PP:=\{1,2,3,\ldots\}.$
The cardinality of a set $A$ will be denoted by $|A|$ and we let ${[n] \choose 2}:=\{S \subseteq [n]: |S|=2 \}.$ For any $n, m \in \ZZ$, we denote $n \equiv m$ if $n \equiv m$ (mod 2). 
For $n \in \PP$ we let,
\[ [n]_{q}:=\frac{1-q^{n}}{1-q}.\] 
It is immediate to see that $[1]_q=[1]_{-q}$ and that for $n\geq 2$
\begin{equation}\label{fattoriale}
[n]_q=\left \{\begin{array}{ll} 
[n]_{-q}+2q+2q^2+\ldots +2q^{n-1} & \mbox{\rm if $n$ is even}, \\
 
[n]_{-q}+2q+2q^2+\ldots +2q^{n-2} & \mbox{\rm if $n$ is odd}. \\
\end{array} \right. 
\end{equation}

\subsection{Symmetric group}

Let $S_{n}$ be the set of all bijections $\sigma :
[n]\rightarrow[n]$.  If $\sigma \in S_{n}$ then we write $\sigma =
\sigma _{1} \ldots \sigma _{n}$ to mean that $\sigma (i) =\sigma
_{i}$, for $i=1, \ldots ,n$. For $\s \in S_n$ and in general, for any sequence $\sigma=(\s_{1},\ldots,\s_{n}) \in \ZZ^{n}$ we
say that $(i,j) \in [n]\times[n]$ is an {\em inversion} of
$\sigma$ if $i<j$ and $\s_{i}>\s_{j}$. We denote the number of inversions of $\s$ by $\inv(\s)$. It is well known that $S_n$ is a Coxeter group respect to the generating set $S:=\{s_i: =(i,i+1) \; : \; i \in [n-1] \}$. The {\em length} of $\s \in S_n$, respect to $S$ is denoted by $\ell(\s)$. It is well known that $\ell(\s)=\inv(\s)$ and the {\em Poincar\'e polynomial} of $W$ is given by  
\[\sum_{\s \in S_n} q^{\ell(\s)}=[1]_{q}[2]_q\cdots [n]_q.\]
We say that $i \in [n-1]$ is a {\em descent} of $\sigma=(\s_{1},\ldots,\s_{n}) \in \ZZ^{n}$ if $\s_{i}>\s_{i+1}$. We denote by
$\Des(\sigma)$ the set of descents and by $\des(\s)$ its cardinality. We also let
\begin{equation}
\label{maj} \maj(\sigma):=\sum_{i \in \Des(\sigma)} i
\end{equation}
and call it the {\em major index} of $\sigma$. 

The definition of Mahonian statistic, come from the following theorem due to MacMahon \cite{MM}. Foata gave a bijective proof of this result in \cite{Fo}.

\begin{thm}[MacMahon] Let $n \in \PP$. Then
\[\sum_{\s \in S_n} q^{\maj(\s)}=\sum_{\s \in S_n} q^{\ell(\s)}.\]
\end{thm}

For any element $w$ of a Coxeter group $W$, the {\em sign} of $w$ is defined to be
\[\sign(w):=(-1)^{\ell(w)}.\]
We prefer to use the notation $(-1)^{\ell(w)}$, instead of the usual $\sign(w)$, in order to avoid confusion between signed and even-signed permutations in the following Section 4. The following is a formula for the signed Mahonian for $S_n$ and is due to Gessel and Simion \cite{W}.
\begin{thm}[Gessel-Simion]\label{GS} Let $n\in \PP$. Then
\[\sum_{\s \in S_n} (-1)^{\ell(\s)}q^{\maj(\s)}=[1]_{q}[2]_{-q}\cdots [n]_{(-1)^{n-1}q}.\]
\end{thm}

\subsection{Hyperoctahedral group}

We denote by $B_{n}$ the group of all bijections $\be$ of the set
$[-n,n]\setminus \{0\}$ onto itself such that
\[\be(-i)=-\be(i)\]
for all $i \in [-n,n]\setminus \{0\}$, with composition as the
group operation. This group is usually known as the group of {\em
signed permutations} on $[n]$, or as the {\em hyperoctahedral
group} of rank $n$.  If
$\be \in B_{n}$ then we write $\be=[\be_{1},\dots,\be_{n}]$ to
mean that $\be(i)=\be_{i}$ for $i=1,\ldots,n$, we call this the
{\em window} notation of $\be$. As set of generators for $B_n$ we
take $S_B:=\{s_1^B,\ldots,s_{n-1}^B,s_0^B\}$, where for $i
\in[n-1]$
\[s_i^B:=[1,\ldots,i-1,i+1,i,i+2,\ldots,n] \;\; {\rm and} \;\; s_0^B:=[-1,2,\ldots,n].\]
It's well known that $(B_n,S_B)$ is a
Coxeter system of type $B$ (see e.g., \cite{BB}). As for $S_n$ we
give an explicit combinatorial description of the length function
$\ell_B$ of $B_n$ with respect to $S_B$. For $\be \in B_n$ we let
$\Neg(\be)   :=  \{i \in [n]: \be_{i}<0\},$
\[\me_{1}(\be) :=  |\Neg(\be)|,\]
and \[\me_{2}(\be)  :=  |\{\{i,j\} \in {[n] \choose
2}:\be_{i}+\be_{j}<0 \}|.\] Note that, if $\be \in B_n$, then it's
not hard to see that
\begin{equation} \label{som}
\me_{1}(\be)+\me_{2}(\be)=-\sum_{i \in \Neg(\be)} \be(i).
\end{equation}
We have the following characterization of the length function (see e.g., \cite{Br}).
\begin{prop} Let $\be \in B_n$. Then
\[\ell_B(\be)=\inv(\be)+\me_1(\be)+\me_2(\be).\]
\end{prop}
The Poincar\'e polynomial of $B_n$ is
\begin{equation}\label{poincareB}
\sum_{\be \in B_{n}} q^{\ell_B(\be)}=[2]_q[4]_q\cdots [2n]_q.
\end{equation}

For any $\be \in B_n$, the {\em flag-major index} of $\be$, here denoted by $\fmaj(\be)$, is defined by 
\begin{equation} \label{defmaj}
\fmaj(\be)=2 \maj(\be)+\me_1(\be),
\end{equation}
where $\maj$ is computed by using the following order on $\ZZ$
\begin{equation}\label{reverseorder} 
-1\prec -2\prec \cdots \prec -n\prec \cdots \prec 0\prec 1\prec 2\prec \cdots \prec n\prec \cdots  
\end{equation}
instead of the usual ordering, $\leq$. 

\noindent For example, if $\be=[2,-5,-3,-1,4] \in B_5$, then $\Des(\be)=\{1,2,3\}$, hence $\maj(\be)=6$ and $\fmaj(\be)=15$. However be aware that $\inv(\be)=3$. 

\smallskip

The $\fmaj$ is a Mahonian statistic on $B_n$ (see \cite[Theorem 2]{AR}).
\begin{thm}[Adin-Roichman]
Let $n \in \PP$. Then
\[ \sum_{\be \in B_{n}} q^{\fmaj(\be)}=\sum_{\be \in B_{n}} q^{\ell_B(\be)}.\]
\end{thm}

The following formula for the signed Mahonian polynomial of type $B$  has been recently discovered by Adin, Gessel, and Roichman \cite{AGR}.

\begin{thm}[Adin-Gessel-Roichman]\label{AGR} Let $n \in \PP$. Then
\[\sum_{\be \in B_n} (-1)^{\ell_B(\be)}q^{\fmaj(\be)}=[2]_{-q}[4]_q\cdots [2n]_{(-1)^nq}.\]
\end{thm}

The group $B_n$ has four one-dimensional characters. We have already shown formulas for the trivial and the sign character. The other two characters are $(-1)^{\me_1(\be)}$ and the sign of $(|\be_1|,\ldots, |\be_n|)$. Their generating functions can be easily obtained from (\ref{poincareB}) and Theorem \ref{AGR}, see \cite{AGR}.
\begin{thm} Let $n \in \PP$. Then
\[\sum_{\be \in B_n} (-1)^{\me_1(\be)} q^{\fmaj(\be)}=[2]_{-q}[4]_{-q}\cdots [2n]_{-q};\]
\[\sum_{\be \in B_n} (-1)^{\ell_B(|\be|)} \; q^{\fmaj(\be)}=[2]_{q}[4]_{-q}\cdots [2n]_{(-1)^{n-1}q}.\]
\end{thm}

\subsection{Even-signed permutation group}

We denote by $D_{n}$ the subgroup of $B_{n}$ consisting of all the
signed permutations having an even number of negative entries in
their window notation, more precisely
\[ D_{n} := \{\g \in B_{n}: \me_{1}(\g)\equiv 0 \}. \]
As a set of generators for $D_n$ we take
$S_D:=\{s_{0}^D,s_{1}^D,\dots,s_{n-1}^D\}$, where for $i \in [n-1]$
\[s_i^D:=s_i^B \;\; {\rm and} \;\; s_{0}^D:=[-2,-1,3,\ldots,n].\]
There is a well
known direct combinatorial way to compute the length for $\g \in
D_{n}$ (see, e.g., \cite[\S 8.2]{BB}). Let $\g \in D_n$. Then
\[ 
\ell_D(\g)=\inv(\g) + \me_{2}(\g).
\]
Note that $\ell_D(\g)=\ell_B(\g)-\me_1(\g)$.
The Poincar\'e polynomial of $D_n$ is
\[\sum_{\g \in D_{n}} q^{\ell_D(\g)}=[2]_q[4]_q\cdots [2n-2]_q[n]_q.\]

For any $\gamma \in D_n$ let 
\[|\gamma|_n:=[\g(1),\ldots,\g(n-1),|\g(n)|].\]
Following \cite{BC}, for $\gamma \in D_n$ we define the \emph{D-major index} by
\[\Dmaj(\gamma):=\fmaj(|\gamma|_n) . \]
We introduce the following subset of $B_n$, 
\[\Delta_n:=\{\g \in B_n : \g(n)>0\}.\] 
The map $\varphi:D_n \longrightarrow  \Delta_n$ defined by $\g \mapsto |\g|_n$ is a bijection. So, by means of this bijection, any function defined on $\Delta_n$ can also be considered as a function defined on $D_n$. In what follows, we work with the subset $\Delta_n$ instead of $D_n$ in order to make some definitions and arguments more natural and transparent. In particular, for $\g \in \Delta_n$, we let 
\begin{equation}\label{ridef}
\Dmaj(\g):=\fmaj(\g) .
\end{equation}

\noindent For example if $ \gamma =[2,-5,-3,-1,4] \in \Delta_5$, then $\varphi(\g)=[2,-5,-3,-1,-4]\in D_5$ and  $\Dmaj(\gamma )=\Dmaj(\varphi(\g))=15.$ 

\smallskip

The statistic $\Dmaj$ is Mahonian on $D_n$ (see \cite[Proposition 4.2]{BC}).
\begin{thm}[Biagioli-Caselli]\label{equiD}
Let $n \in \PP$. Then
\[\sum _{\gamma \in \Delta_{n}}q^{\Dmaj(\gamma )}=\sum _{\gamma \in D_{n}}q^{\ell_D(\gamma )}.\]
\end{thm}

\section{A sign-reversing involution on $B_n$}

In this section we define a natural involution on $B_n$ and derive
some of its properties. In particular this allows an easy proof of the
Adin-Gessel-Roichman formula for the signed-Mahonian for $B_{2n}$. We will limit our
discussion to the involution for $B_{2n}$. The odd case is almost identical.

\smallskip

Let $\iota: B_{2n} \rightarrow B_{2n}$ be the map defined by
\begin{equation} \label{definvolution}
\be \mapsto s_{2i-1} \cdot \be,
\end{equation}
where $i \in [n]$ is the smallest integer such that $2i-1$ and $2i$
have opposite signs or are not in adjacent positions in the windows
notation of $\be$. It no such $i$ exists let $\iota(\be):=\beta$.

\noindent For example, $\be=[-3,-4,1,2,-6,-5]$ is a fixed point, and
$\g=[2,6,5,-4,-3,1]$ is such that $\iota(\g)=s_1\g=[1,6,5,-4,-3,2].$

It is clear that when $\be \in B_{2n}$ is not a fixed point, the involution $\iota$ reverses the sign of $\be$. However, it preserves
the descent set $\Des(\be)$, the number of negative entries $\me_1(\be)$, and hence the flag-major index $\fmaj(\be)$. Namely, 
\begin{equation}\label{signrev}
\ell_B(\be)\not \equiv \ell_B(\iota(\be)) \;\; {\rm and} \;\; \fmaj(\be)=\fmaj(\iota(\be)).
\end{equation}

The following lemma will be fundamental in the proof of our main
result.

\begin{lem}\label{quarto}
Let $n \in \PP$. Then
\begin{equation}\label{quartoeq}
\sum_{\genfrac{}{}{0pt}{}{\be \in B_{2n}}{\fmaj(\be)\equiv 1}}(-1)^{\ell_B(\be)}q^{\fmaj(\be)}=0.
\end{equation}
\end{lem}
\begin{proof} 
From (\ref{signrev}), all the terms in the RHS of (\ref{quartoeq}) cancel
except for the terms corresponding to the fixed points. Now let $\be
\in B_{2n}$ such that $\fmaj(\be)\equiv 1$. This implies that in the window
notation of $\be$ there is an odd number of both negative and positive
entries. Hence there exists at least one pair $2i-1$, $2i$ in $\be$
with opposite signs. It follows that $\be$ is not a fixed point and
this concludes the proof.
\end{proof}

Now, let consider the set of fixed point of $\iota$. Such elements in $B_{2n}$
correspond bijectively to signed permutations in $B_n$
with some entries ``barred'' according with the following rule:
\begin{itemize}
\item[] each pair of adjacent entries of type $\pm (2i-1)$, $\pm \: 2i$ in $\be$ is
replaced by $\pm \: i$;
\item[] each pair of adjacent entries of type $\pm \: 2i$, $\pm (2i-1)$ in $\be$ is
replaced by $\pm \: \bar{i}$.
\end{itemize}

\noindent We denote by $\mathcal{B}_n$ the set of all the {\em barred signed permutations} of $[\pm n]$. Let $\Des(\bar{\be})$ and $\maj(\bar{\be})$ be defined without considering the bars and let $\SO(\bar{\be})$ be the set of the positions of the barred entries.

\noindent For example, let $\be=[-3,-4,1,2,-6,-5]\in B_6$ be a fixed point. Then $\bar{\be}=[-2,1,-\bar{3}] \in \mathcal{B}_3$ is the corresponding barred signed permutation, with $\fmaj(\bar{\be})=5$ and $\SO(\bar{\be})=\{3\}$.

Since to compute the descent set we consider the reverse ordering (\ref{reverseorder}) on the negative integers, it follows that
\[\Des(\be)=\{2i \; : \; i \in \Des(\bar{\be})\} \cup \{2i-1 \; : \; i \in \SO(\bar{\be})\}.\]
Hence
\begin{equation}\label{majeq}
\maj(\be)= 2 \maj(\bar{\be}) + \sum_{i \in \SO(\bar{\be})} (2i -1).
\end{equation} 
Differently, to compute the inversions of $\be$ we use the natural order $\leq$. Hence 
\begin{eqnarray}\label{lunghezza}
\inv(\be)& = & 4\inv(\bar{\be})+|\SO(\bar{\be})^+|+(\; \me_1(\bar{\be})-|\SO(\bar{\be})^-| \;) \nonumber \\
& \equiv & \me_1(\bar{\be})+|\SO(\bar{\be})|,  
\end{eqnarray}
where $|\SO(\bar{\be})^{\pm}|$ denote the number of positive and negative barred entries in $\bar{\be}$, respectively. Moreover, it is easy to see that  
\begin{eqnarray}\label{lunghezza1}
\me_1(\be)+\me_2(\be)=-\sum_{i \in \Neg(\bar{\be})}(4\bar{\be}_i+1).
\end{eqnarray}

The following theorem holds.

\begin{thm}\label{ricorsionepari}
Let $n \in \PP$. Then
\[ \sum_{\be \in B_{2n}}(-1)^{\ell_B(\be)}q^{\fmaj(\be)} =\prod_{i=1}^n(1-q^{4i-2})\sum_{\be \in B_n}q^{2\fmaj(\be)} .\]
\end{thm}
\begin{proof}
Since the involution $\iota$ is sign-reversing and $\fmaj$ preserving, to compute the LHS is enough to perform the sum over the set of fixed points of $\iota$. So let $\be \in B_{2n}$ a fixed point and $\bar{\be}$ the corresponding barred signed permutation in $\mathcal{B}_n$. From (\ref{majeq}), (\ref{lunghezza}), (\ref{lunghezza1}) and $\me_1(\be)=2\me_1(\bar{\be})$, it follows 
\[
\fmaj(\be)=4 \maj(\bar{\be}) + \sum_{i \in \SO(\bar{\be})} (4i -2) + 2 \me_1(\bar{\be})\]
and 
\[\ell_B(\be)\equiv\me_1(\bar{\be})+|\SO(\bar{\be})|+\sum_{i \in \Neg(\bar{\be})}(4\bar{\be}_i+1) \equiv |\SO(\bar{\be})|.\]
Hence 
\begin{eqnarray*}
\sum_{\be \in B_{2n}}(-1)^{\ell_B(\be)}q^{\fmaj(\be)} & = & \sum_{\bar{\be}\in \mathcal{B}_n}(-1)^{|\SO(\bar{\be})|}q^{\sum_{i \in \SO(\bar{\be})}(4i-2)}q^{2\fmaj(\bar{\be})}\\
& = & \prod_{i=1}^n (1 - q^{4i -2})\sum_{\be \in B_n}q^{2\fmaj(\be)}.
\end{eqnarray*}
\end{proof}

The case even of Theorem \ref{AGR} follows directly by this and (\ref{poincareB}).

\begin{cor} Let $n \in \PP$. Then
\[\sum_{\be \in B_{2n}}(-1)^{\ell_B(\be)}q^{\fmaj(\be)}=[2]_{-q}[4]_q \ldots [4n]_q.\]
\end{cor} 

In the case of $B_{2n+1}$, the fixed points of the involution $\iota$ are the signed permutations such that, for every $i\in [n]$ the entries $2i-1$, $2i$ have same sign and are in adjacent positions in the window notation of $\be$. This forces the entry $\pm (2n+1)$ to be in an odd-position. It is possible to obtain the Adin-Gessel-Roichman formula also in this case, again by summing over the set of fixed points just described. Since this procedure it is very technical and not easier than the original proof it is not worth presenting here.

\section{A Signed Mahonian for $D_n$}

In this section we provide a factorial-type formula for the signed Mahonian polynomial for the Weyl group $D_n$.

For any $\be=[\be_1,\ldots,\be_n]$ let $-\be:=[-\be_1,\ldots,-\be_n]$. It is easy to see that the following equalities hold:
\[\inv(-\be) = {n \choose 2} - \inv(\be), \\
\me_1(-\be) = n - \me_1(\be), \;\; {\rm and} \;\; \me_2(-\be) = {n \choose
2} - \me_2(\be).
\]
It follows
\begin{equation}\label{menobeta}
\ell_B(-\be)\equiv \ell_B(\be)+n.
\end{equation}

A proof of the following lemma can be found in \cite[Corollary 3.13]{BC}.
\begin{lem}\label{menofmaj} Let $\g \in \Delta_n$. Then
\[\fmaj(-\g)=\fmaj(\g)+n.\]
\end{lem}

\begin{prop}\label{zero}
Let $n \in \PP$. Then
\[\sum_{\be \in B_n} (-1)^{\ell_B(\be)}q^{\fmaj(\be)}=
\sum_{\g \in \Delta_n}
(-1)^{\ell_B(\g)}q^{\Dmaj(\g)}(1+(-q)^n).\]
\end{prop}
\begin{proof}
From Lemma \ref{menofmaj}, (\ref{menobeta}) and (\ref{ridef}) we have
\begin{eqnarray*}
\sum_{\be \in B_n} (-1)^{\ell_B(\be)}q^{\fmaj(\be)}& = & \sum_{\g
\in \Delta_n} (-1)^{\ell_B(\g)}q^{\fmaj(\g)}+
(-1)^{\ell_B(-\g)}q^{\fmaj(-\g)} \\
& = & \sum_{\g \in \Delta_n} (-1)^{\ell_B(\g)}q^{\fmaj(\g)}+
(-1)^{\ell_B(\g)+n}q^{\fmaj(\g)+n} \\
& = & \sum_{\g \in \Delta_n}
(-1)^{\ell_B(\g)}q^{\fmaj(\g)}(1+(-q)^n)\\
& = &\sum_{\g \in \Delta_n}
(-1)^{\ell_B(\g)}q^{\Dmaj(\g)}(1+(-q)^n).
\end{eqnarray*}
\end{proof}
The following are immediate consequences of Theorem \ref{AGR}.
\begin{cor}\label{pari} Let $n\in \PP$ be even. Then
\[\sum_{\g \in \Delta_n}(-1)^{\ell_B(\g)}q^{\Dmaj(\g)}=[2]_{-q}[4]_{q}\cdots[2n-2]_{-q}[n]_{q}.\]
\end{cor}
\begin{cor}\label{dispari} Let $n\in \PP$ be odd. Then
\[\sum_{\g \in \Delta_n}(-1)^{\ell_B(\g)}q^{\Dmaj(\g)}=[2]_{-q}[4]_{q}\cdots[2n-2]_{q}[n]_{-q}.\]
\end{cor}

We denote by $\Delta_n^0$ and $\Delta_n^1$ the subsets of all $\g \in \Delta_n$ such that $\Dmaj(\g)\equiv 0$ and  $\Dmaj(\g)\equiv 1$, respectively. The subsets $D_n^0$ and $D_n^1$ are defined in a similar way. 

\begin{lem}\label{primo}
Let $n \in \PP$. Then
\[\sum_{\g \in \Delta_n^0}
(-1)^{\ell_B(\g)}q^{\Dmaj(\g)}=\sum_{\g \in D_n^0}
(-1)^{\ell_D(\g)}q^{\Dmaj(\g)}.\]
\end{lem}
\begin{proof}
Every signed permutation $\g \in \Delta_n^0$ is such that $\me_1(\g)\equiv 0$.
Hence $\g \in D_n^0$ and  
\[\ell_B(\g)=\inv(\g)+\me_1(\g)+\me_2(\g)\equiv
\inv(\g)+\me_2(\g)=\ell(\g).\]
\end{proof}

\begin{lem}\label{secondo}
Let $n \in \PP$. Then \begin{equation}\label{toprove} \sum_{\be
\in\Delta_n^1} (-1)^{\ell_B(\be)}q^{\Dmaj(\be)}=-\sum_{\g \in
D_n^1} (-1)^{\ell_D(\g)}q^{\Dmaj(\g)}.
\end{equation}
\end{lem}
\begin{proof} Let $\be:=[\be_1, \ldots,\be_{n-1},k] \in \Delta_n^1$. This implies
$\me_1(\be)\equiv 1$ and so $\be \in B_n
\setminus D_n$. Let $\g:=\varphi^{-1}(\be)$, i.e., $\g:=[\be_1, \ldots,\be_{n-1},-k] \in D_n^1$. We have
\[\me_2(\g)=\me_2(\be)+(k-1) \;\ {\rm and} \;\;
\inv(\g)=\inv(\be)+(k-1).\] 
It follows that
\[\ell_B(\be)=\inv(\be)+\me_1(\be)+\me_2(\be)\not\equiv
\inv(\be)+2(k-1)+\me_2(\be)=\ell_D(\g).\] 
\end{proof}

\begin{lem}\label{terzo} Let $n \in \PP$ even. Then
\begin{equation}\label{eqterzo}
\sum_{\be \in \Delta_n^1} (-1)^{\ell_B(\be)}q^{\Dmaj(\be)}=0
\end{equation}
\end{lem}
\begin{proof} Let consider the restriction of the involution $\iota:B_n \rightarrow B_n$ defined in (\ref{definvolution}) to $\Delta_n^1$. It is easy to see that none of the elements of $\Delta_n^1$ is a fixed point for $\iota$, and that $\Dmaj(\be)=\Dmaj(\iota({\be}))$. Hence all the terms in the RHS of (\ref{eqterzo}) cancel and the result follows.
\end{proof}

Now, we are ready to show the main result of this section.

\begin{thm}[Signed Mahonian of type $D$] Let $n \in \PP$. Then
\[ \sum_{\g \in D_n} (-1)^{\ell_D(\g)}q^{\Dmaj(\g)}= \left \{\begin{array}{ll}
[2]_{-q}[4]_{q}\cdots[2n-2]_{-q}[n]_{q} & \mbox{\rm if $n$ is even,} \\
 
[2]_{-q}[4]_{q}\cdots[2n-2]_{q} \; [n]_{q} & \mbox{\rm if $n$ is odd.}   
\end{array} \right. \]
\end{thm}
\begin{proof} If $n$ is even, from Lemmas \ref{primo} and \ref{secondo}, and Lemma \ref{terzo}, it follows
\[\sum_{\g \in D_n}(-1)^{\ell_D(\g)}q^{\Dmaj(\g)}=\sum_{\g\in \Delta_n}(-1)^{\ell_B(\g)}q^{\Dmaj(\g)}.\]
Hence the result follows by Corollary \ref{pari}.

\smallskip

If $n$ is odd, from Corollary \ref{dispari}, we have
\[\sum_{\be \in \Delta_n}(-1)^{\ell_B(\be)}q^{\Dmaj(\be)}=[2]_{-q}[4]_{q}\cdots[2n-2]_q[n]_{-q}.\]
By Theorem \ref{AGR}, this implies 
\begin{equation}\label{RHS}
\sum_{\be \in \Delta_n}(-1)^{\ell_B(\be)}q^{\Dmaj(\be)}=\sum_{\be \in
B_{n-1}}(-1)^{\ell_B(\be)}q^{\fmaj(\be)}\cdot[n]_{-q}.
\end{equation}
By Corollary \ref{quarto}, the first factor
in the RHS of (\ref{RHS}) has only even powers, hence
\begin{equation}\label{doppio}
\sum_{\be \in \Delta^1_n}(-1)^{\ell_B(\be)}q^{\Dmaj(\be)}=\sum_{\be \in
B_{n-1}}(-1)^{\ell_B(\be)}q^{\fmaj(\be)}\cdot(-q-q^3-\ldots
-q^{n-2}).
\end{equation}
Again from Lemmas \ref{primo} and \ref{secondo} and (\ref{RHS}), it follows
\begin{eqnarray*}
\sum_{\g \in D_n}(-1)^{\ell_D(\g)}q^{\Dmaj(\g)} & = & \sum_{\be
\in B_{n-1}}(-1)^{\ell_B(\be)}q^{\fmaj(\be)}\cdot [n]_{-q}\\
& - & 2 \cdot \sum_{\be \in \Delta_n^1}(-1)^{\ell_B(\be)}q^{\Dmaj(\be)}
\end{eqnarray*}
Now by (\ref{doppio}), (\ref{fattoriale}) and Theorem \ref{AGR} the RHS is equal to
\begin{eqnarray*}
& = & \sum_{\be \in B_{n-1}}(-1)^{\ell_B(\be)}q^{\fmaj(\be)}\cdot
([n]_{-q} + 2(q+q^3+\ldots +q^{n-2}))\\
& = & \sum_{\be
\in B_{n-1}}(-1)^{\ell_B(\be)}q^{\fmaj(\be)}\cdot [n]_{q}\\
& = & [2]_{-q}[4]_{q}\cdots[2n-2]_{q} [n]_{q}.
\end{eqnarray*}
\end{proof}

\bigskip

\noindent LaCIM - Universit\'e du Qu\'ebec \`a Montr\'eal, \\
Case Postale 8888, succursale Centre-ville, Montr\'eal, Qu\'ebec, Canada H3C 3P8
\smallskip

\noindent {\em E-mail}: {\tt biagioli@math.uqam.ca}

\noindent {\em URL}: {\tt http://www.lacim.uqam.ca/\~{\rule{1pt}{0pt}}biagioli}


\begin{thebibliography}{xx}

\bibitem{AR}
R. M. Adin and Y. Roichman, {\em The Flag Major Index and group actions on polynomial rings}, Europ. J. Combinatorics, {\bf 22} (2001), 431-446.

\bibitem{ABR} R. M. Adin, F. Brenti and Y. Roichman, \emph{Descent representations and multivariate statistics}, Trans. Amer. Math. Soc., to appear.

\bibitem{ABR1} R. M. Adin, F. Brenti and Y. Roichman, \emph{Descent numbers and major indices for the hyperoctahedral group}, Adv. Appl. Math., {\bf 27} (2001), 210-224. 

\bibitem{AGR}
R. M. Adin, I. Gessel and Y. Roichman, {\em Signed Mahonians}, preprint 2004, arXiv:math.CO/0402208.

\bibitem{BC} R. Biagioli and F. Caselli, \emph{Invariant algebras and major indices
for classical Weyl groups,} Proc. London Math. Soc., {\bf 88} (2004), 603-631.

\bibitem{BC2} R. Biagioli and F. Caselli, \emph{A descent basis for the coinvariant algebra of type $D$,} J. Algebra, {\bf 275} (2004), 517-539.

\bibitem{BB}
A. Bj\"{o}rner and F. Brenti, {\em Combinatorics of Coxeter groups},
to appear on Grad. Texts in Math., Springer-Verlag, Berlin, 2001

\bibitem{Br}
F. Brenti, {\em q-Eulerian polynomials arising from Coxeter groups}, Europ. J. Combinarorics, {\bf 15} (1994), 417-441.

\bibitem{Ca}
L. Carlitz, {\em A combinatorial property of q-Eulerian numbers}, Amer. Math. Montly, {\bf 82} (1975), 51-54.

\bibitem{CG}
C. O. Chow and I. Gessel, {\em On descent numbers and major indices for the hyperoctahedral group}, preprint 2003.

\bibitem{DF}
J. D\'esarm\'enien and D. Foata, {\em The signed Eulerian numbers}, Discrete Math, {\bf 99} (1992), 49-58.

\bibitem{Fo}
D. Foata, {\em On the Netto inversion number of a sequence}, Proc. Amer. Math. Soc., {\bf 19} (1968), 236-240.

\bibitem{GG}
A. Garsia and I. Gessel, {\em Permutation statistics and partitions}, Adv. Math., {\bf 31} (1979), 288-305.

\bibitem{Hum}
J. E. Humphreys, {\em Reflection groups and Coxeter groups},
Cambridge Stud. Adv. Math., n. 29, Cambridge Univ. Press,
Cambridge,  1990.

\bibitem{MM}
P. A. MacMahon, {\em Combinatory analysis}, vol. 1, Cambridge Univ. Press, London, 1915.

\bibitem{Re}
V. Reiner, {\em Descent and one-dimensional characters for classical Weyl groups}, Discrete Math., {\bf 140} (1995), 129-140.

\bibitem{Re1}
V. Reiner, {\em Signed permutation statistics}, Europ. J. Combinatorics, {\bf 14} (1993), 553-567.

\bibitem{StaEC1}
 R. P. Stanley,
{\em Enumerative combinatorics}, vol. 1, Wadsworth and
Brooks/Cole, Monterey, CA, 1986.

\bibitem{W}
M. L. Wachs, {\em An involution for signed Eulerian numbers}, Discrete Math, {\bf 99} (1992), 59-62.

\end{thebibliography}
\end{document}